\documentclass{amsart2000}

\usepackage[backref=section]{hyperref}

\author{Gergely Harcos}
\title{Uniform approximate functional equation for principal $L$-functions}
\address{Department of Mathematics, Princeton University, Fine Hall,
Washington Road, Princeton, NJ 08544, USA}
\email{gharcos@math.princeton.edu}

\subjclass[2000]{11M41, 11F66}

\newtheorem{Th}{Theorem}
\newtheorem{Cor}{Corollary}
\newtheorem{Lemma}{Lemma}

\newcommand{\al}{\alpha}
\newcommand{\de}{\delta}
\newcommand{\ep}{\epsilon}
\newcommand{\GG}{\Gamma}
\newcommand{\ka}{\kappa}
\newcommand{\la}{\lambda}
\newcommand{\si}{\sigma}

\newcommand{\CC}{\mathbb{C}}

\newcommand{\RR}{\mathbb{R}}

\begin{document}

\begin{abstract}
We prove an approximate functional equation for the central value
of the $L$-series attached to
an irreducible cuspidal automorphic representation $\pi$ of $GL_m$ over
a number field with unitary central character and contragradient
representation $\tilde\pi$.
The approximation involves a smooth truncation of the Dirichlet series
$L(s,\pi)$ and $L(s,\tilde\pi)$ after about $\sqrt{C}$ terms,
$C=C(\pi)=C(\tilde\pi)$ being
the analytic conductor introduced by Iwaniec and Sarnak \cite{IS}.
We investigate the decay rate of the cutoff function and its derivatives
(Theorem~\ref{Th1}). We also see
that the truncation can be made uniformly explicit at the
cost of an error term (Theorem~\ref{Th2}).
Straightforward extensions of these results exist for
products of central values.
We hope that these formulae will help further understanding of the
central values of principal $L$-functions, such as
finding good bounds on their various power means, or
establishing subconvexity or nonvanishing results in certain families.
\end{abstract}

\maketitle

\section{Introduction}

In their discussion \cite{IS} of families of $L$-functions and the
corresponding (sub)con\-vex estimates Iwaniec and Sarnak
introduced the analytic conductor $C=C(\pi)$ of a cusp form $\pi$
on $GL_m$ over a number field $F$. As they pointed out, the
Phragm\'en--Lindel\"of principle implies that the central value
$L(1/2,\pi)$ is at most $C^{1/4+\ep}$, an estimate referred to as
the convexity or trivial bound. The central value contains
important arithmetic information, so it is often very useful (in
fact crucial) to replace the exponent $1/4$ by any smaller value
(note that the generalized Lindel\"of hypothesis which in turn is
implied by the generalized Riemann hypothesis asserts that any
positive exponent is permissible). It is sometimes possible to
realize this improvement by placing the $L$-function into a family
or a fake family and averaging some moment of the corresponding
central values with well-chosen weights (called amplifiers). A
crucial ingredient in such an argument is to express the central
values $L(1/2,\pi)$ in the family by an ``approximate functional
equation'', i.e. a sum of two Dirichlet series which essentially
have $\sqrt{C}$ terms. For details we refer the reader to
\cite{IS}.

In the present note we try to perform the calculation for the
entire family of $\pi$'s on $GL_m$ over $F$ ($m$ and $F$ are
fixed). First we obtain an exact representation of the central
value with uniform decay properties (Theorem~\ref{Th1}). This
formula is most useful for families whose Archimedean parameters
remain bounded. The second representation (Theorem~\ref{Th2}),
inspired by a recent result of Ivi\'c \cite{Iv}, has a more
explicit main term at the cost of an error term. This formula
works best in families where the Archimedean parameters
simultaneously grow large.

The proofs are based on standard Mellin transform techniques
combined with recent progress on the Ramanujan--Selberg
conjectures achieved by Luo, Rudnick and Sarnak \cite{LRS}. More
precisely, we use the bounds of \cite{LRS} at the Archimedean
places to see that the Mellin integrals can be evaluated in a
sufficiently large half plane, while the bounds of \cite{LRS} at
the non-Archimedean places enter through the work of Molteni
\cite{M} by providing the necessary estimates for the Dirichlet
coefficients of the $L$-functions. A variant of the method yields
similar formulae for products of central values (e.g. for higher
moments).

\section{Statement of results}

Let $F$ be a number field of degree $d$ and $\pi=\otimes_v\pi_v$
be an irreducible cuspidal automorphic representation of $GL_m$ over $F$
with unitary central character and contragradient representation $\tilde\pi$.
The corresponding
$L$-functions are defined for $\Re s>1$ by absolutely convergent
Dirichlet series as
\begin{equation}\label{eq3}
L(s,\pi)=\sum_{n=1}^\infty \frac{a_n}{n^s}
\qquad\text{and}\qquad
L(s,\tilde\pi)=\sum_{n=1}^\infty \frac{\overline{a}_n}{n^s},
\end{equation}
and these are connected by a functional equation of the form
\begin{equation}\label{eq2}
N^\frac{s}{2}L(s,\pi_\infty)L(s,\pi)= \ka
N^\frac{1-s}{2}L(1-s,\tilde\pi_\infty)L(1-s,\tilde\pi).
\end{equation}
Here $N$ is the conductor (a positive integer), $\ka$ is the root
number (of modulus 1) and
\begin{equation}\label{eq7}
L(s,\pi_\infty)=\prod_{j=1}^{md}\pi^{-\frac{s}{2}}\GG\left(\frac{s+\mu_j}{2}\right),
\quad L(s,\tilde\pi_\infty)=\prod_{j=1}^{md}\pi^{-\frac{s}{2}}
\GG\left(\frac{s+\overline{\mu}_j}{2}\right)
\end{equation}
are the products of the $L$-functions of $\pi_v$ and
$\tilde\pi_v$, respectively, at the Archimedean places $v$. We
define the analytic conductor of $\pi$ (at $s=1/2$) as
\begin{equation}\label{eq23}
C=\frac{N}{\pi^{md}}\prod_{j=1}^{md}\left|\frac{1}{4}+\frac{\mu_j}{2}\right|.
\end{equation}
(Hopefully no confusion arises from the fact that $\pi$ denotes both a
representation and a real constant.
The meaning should be clear from the context.)
We have the following uniform approximate functional equation
expressing $L(1/2,\pi)$ as a sum of two Dirichlet series.

\begin{Th}\label{Th1} There is a smooth function $f:(0,\infty)\to\CC$ and a complex
number $\la$ of modulus 1 depending only on the Archimedean
parameters $\mu_j$ $(j=1,\dots,md)$ such that
\begin{equation}\label{eq8}
L\left(\frac{1}{2},\pi\right)=
\sum_{n=1}^\infty\frac{a_n}{\sqrt{n}}f\left(\frac{n}{\sqrt{C}}\right)+
\ka\la\sum_{n=1}^\infty\frac{\overline{a}_n}{\sqrt{n}}
\overline{f}\left(\frac{n}{\sqrt{C}}\right).
\end{equation}
The function $f$ and its partial derivatives $f^{(k)}$
$(k=1,2,.\dots)$ satisfy the following uniform growth estimates at
$0$ and infinity:
\begin{equation}\label{eq5}f(x)=
\begin{cases}1+O_{\si}(x^{\si}),&\quad
0<\si<1/(m^2+1);\\
O_{\si}(x^{-\si}),&\quad \si>0;
\end{cases}
\end{equation}
\begin{equation}\label{eq12}
f^{(k)}(x)=O_{\si,k}(x^{-\si}),\quad \si>k-\frac{1}{m^2+1}.
\end{equation}
The implied constants depend only on $\si$, $k$, $m$ and $d$.
\end{Th}

\noindent{\bf Remark 1.}
The range $0<\si<1/(m^2+1)$ in (\ref{eq5}) can be widened
to $0<\si<1/2$ for all representations $\pi$ which are tempered at
$\infty$, i.e. conjecturally for all $\pi$. Similarly, upon the
Ramanujan--Selberg conjecture the range of $\si$ in (\ref{eq12})
can be extended to $\si>k-1/2$.

\medskip

The following corollaries are simple consequences of
Theorem~\ref{Th1} combined with an average form of the (finite)
Ramanujan conjecture recently obtained by Molteni (Theorem 4 of
\cite{M}):
\begin{equation}\label{eq17}
\sum_{n\leq x}|a_n|=O_{\ep,m,d}(x^{1+\ep}C^\ep).
\end{equation}

\begin{Cor}\label{Cor1} For any positive numbers $\ep$ and $A$,
\[L\left(\frac{1}{2},\pi\right)=
\sum_{n\leq
C^{1/2+\ep}}\frac{a_n}{\sqrt{n}}f\left(\frac{n}{\sqrt{C}}\right)+
\ka\la\sum_{n\leq C^{1/2+\ep}}\frac{\overline{a}_n}{\sqrt{n}}
\overline{f}\left(\frac{n}{\sqrt{C}}\right)+O_{\ep,A}(C^{-A}).\]
The implied constant depends only on $\ep$, $A$, $m$ and $d$.
\end{Cor}

\begin{Cor}\label{Cor2} For any $\ep>0$, there is a uniform convexity bound
\[L(1/2,\pi)\ll_{\ep}C^{1/4+\ep}.\]
The implied constant depends only on $\ep$, $m$ and $d$.
\end{Cor}

In a family of representations $\pi$ it is often desirable to see
that the weight functions $f$ do not vary too much. In fact,
assuming that the Archimedean parameters are not too small one can
replace $f$ by an explicit function $g$ (independent of $\pi$) and
derive an approximate functional equation with a nontrivial error
term, i.e., an error substantially smaller than the convexity
bound furnished by the above corollary. To state the result we
introduce
\begin{equation}\label{eq19}
\eta=\min_{j=1,\dots,md}\left|\frac{1}{4}+\frac{\mu_j}{2}\right|.
\end{equation}

\begin{Th}\label{Th2} Let $g:(0,\infty)\to\RR$ be a smooth function
with functional equation $g(x)+g(1/x)=1$ and derivatives decaying
faster than any negative power of $x$ as $x\to\infty$. Then, for
any $\ep>0$,
\[L\left(\frac{1}{2},\pi\right)=
\sum_{n=1}^\infty\frac{a_n}{\sqrt{n}}g\left(\frac{n}{\sqrt{C}}\right)+
\ka\la\sum_{n=1}^\infty\frac{\overline{a}_n}{\sqrt{n}}
g\left(\frac{n}{\sqrt{C}}\right)+O_{\ep,g}(\eta^{-1}C^{1/4+\ep}),\]
where $\la$ (of modulus 1) is given by (\ref{eq16}) and the
implied constant depends only on $\ep$, $g$, $m$ and $d$.
\end{Th}

\noindent{\bf Remark 2.} This should be compared with the main
result of Ivi\'c \cite{Iv} (in the light of the next remark). It
is apparent that the formula is really of value when the family
under consideration satisfies $\eta\gg C^{\de}$ with some fixed
$\de>0$.

\medskip

\noindent{\bf Remark 3.} We obtain similar expressions for any
value $L(1/2+it,\pi)$ on the critical line by twisting $\pi$ with
the one-dimensional representation $|\det|^{it}$. Under the twist
the conductor remains $N$, the root number becomes $\ka N^{-it}$,
and the Archimedean parameters change to $\mu_j+it$
$(j=1,\dots,md)$. Accordingly, the analytic conductor $C$ needs to
be adjusted as well. Also, similar formulae hold when $\pi$ is not
irreducible but a tensor product of finitely many irreducible
representations. This observation might be useful for studying
higher moments of the central values in families.

\medskip

\section{Proof of Theorem~\ref{Th1}}

By a result of Luo, Rudnick and Sarnak (Theorem 1 in \cite{LRS}),
\begin{equation}\label{eq21}
\Re\mu_j\geq\frac{1}{m^2+1}-\frac{1}{2},\quad j=1,\dots,md.
\end{equation}
(The Ramanujan--Selberg Conjecture asserts that $\pi_\infty$
is tempered upon which we could replace the right hand side by 0.)
Therefore the function
\begin{equation}\label{eq11}
F(s,\pi_\infty)=\left\{N^s\frac
{L(1/2+s,\pi_\infty)L(1/2,\tilde\pi_\infty)}
{L(1/2-s,\tilde\pi_\infty)L(1/2,\pi_\infty)}
\right\}^{1/2}
\end{equation}
is holomorphic in the half plane $\Re s>-1/(m^2+1)$.
With this notation
we can rewrite the functional equation (\ref{eq2}) as
\begin{equation}\label{eq1}
F(s,\pi_\infty)L(1/2+s,\pi)=\ka\la
F(-s,\tilde\pi_\infty)L(1/2-s,\tilde\pi),
\end{equation}
where
\begin{equation}\label{eq16}
\la=\frac{L(1/2,\tilde\pi_\infty)}{L(1/2,\pi_\infty)}.\end{equation}
Note that $F(0,\pi_\infty)=1$ and
\begin{equation}\label{eq6}
\overline{F}(s,\pi_\infty)=F(\overline{s},\tilde\pi_\infty)
\end{equation}
follow from definitions (\ref{eq11}) and (\ref{eq7}). Similarly,
$\la$ is of modulus 1.

We also fix an entire function $H(s)$ which satisfies
the growth estimate
\begin{equation}\label{eq20}
H(s)\ll_{\si,A}(1+|s|)^{-A},\quad\Re s=\si;
\end{equation}
on vertical lines. In addition, we shall assume that $H(0)=1$ and that
$H(s)$ is symmetric with respect to both axes:
\begin{equation}\label{eq9}
H(s)=H(-s)=\overline{H}(\overline{s}).
\end{equation}
Such a function can be obtained as the Mellin transform of a
smooth function $h:(0,\infty)\to\RR$ with total mass 1 with
respect to the measure $dx/x$, functional equation $h(1/x)=h(x)$
and derivatives decaying faster than any negative power of $x$ as
$x\to\infty$:
\[H(s)=\int_0^\infty h(x)x^s\frac{dx}{x}.\]

Using those two auxiliary functions and taking an arbitrary
$0<\si<1/(m^2+1)$ we can express the central value $L(1/2,\pi)$
via the residue theorem as
\[\begin{split}
L(1/2,\pi)
&=\frac{1}{2\pi i}\int_{(\si)}L(1/2+s,\pi)F(s,\pi_\infty)H(s)\frac{ds}{s}\\\\
&-\frac{1}{2\pi i}\int_{(-\si)}L(1/2+s,\pi)F(s,\pi_\infty)H(s)\frac{ds}{s}.
\end{split}\]
Here we combined inequality (\ref{eq20}), Lemma~\ref{Lemma2}
below, and the fact that $L(1/2+s,\pi)$ grows moderately on the
lines $\Re s=\pm\si$. The last property follows from the
Phragm\'en--Lindel\"of principle.

Here we used that the $L$-functions are bounded on the relevant
lines as well as Lemma~\ref{Lemma2} below which shows that
$F(s,\pi_\infty)$ grows moderately on $\Re s=\pm\si$. Applying a
change of variable $s\mapsto -s$ in the second integral we get, by
the functional equations (\ref{eq1}) and (\ref{eq9}),
\[\begin{split}
L(1/2,\pi)
&=\frac{1}{2\pi i}\int_{(\si)}L(1/2+s,\pi)F(s,\pi_\infty)H(s)\frac{ds}{s}\\\\
&+\frac{\ka\la}{2\pi
i}\int_{(\si)}L(1/2+s,\tilde\pi)F(s,\tilde\pi_\infty)H(s)\frac{ds}{s}.
\end{split}\]
By another change of variable $s\mapsto\overline{s}$ in the second
integral we observe, using (\ref{eq3}), (\ref{eq6}) and
(\ref{eq9}), that this integral is the complex conjugate of the
first one. Therefore we obtain the representation
(\ref{eq8}) of Theorem~\ref{Th1} by defining
\begin{equation}\label{eq25}
f\left(\frac{x}{\sqrt{C}}\right)= \frac{1}{2\pi
i}\int_{(\si)}x^{-s}F(s,\pi_\infty)H(s)\frac{ds}{s}.
\end{equation}
Then, for any nonnegative integer $k$ we also have
\begin{equation}\label{eq13}
f^{(k)}(x)=\frac{(-1)^k}{2\pi i}
\int_{(\si)}x^{-s-k}C^{-s/2}F(s,\pi_\infty)H(s)s(s+1)\dots(s+k-1)
\frac{ds}{s}.
\end{equation}

When $k=0$, the integrand in this expression is holomorphic for
$\Re s>-1/(m^2+1)$ with the exception of a simple pole at $s=0$
with residue 1. So in this case we are free to move the line of
integration to any nonzero $\si>-1/(m^2+1)$ but negative $\si$'s
will pick up an additional value 1 from the pole at $s=0$. When
$k>0$, the integrand is holomorphic in the entire half plane $\Re
s>-1/(m^2+1)$, so the line of integration can be shifted to any
$\si>-1/(m^2+1)$ without changing the value of the integral.
Henceforth, by (\ref{eq20}) and (\ref{eq13}), the truth of
inequalities (\ref{eq5}) and (\ref{eq12}) are reduced to the
following:

\begin{Lemma}\label{Lemma2} For any $\si>-1/(m^2+1)$ we have the
uniform bound
\begin{equation}\label{eq22}
C^{-s/2}F(s,\pi_\infty)\ll_{\si}(1+|s|)^{md\si/2},\quad\Re s=\si.
\end{equation}
The implied constant depends only on $\si$, $m$ and $d$.
\end{Lemma}

We start with the following simple estimate.

\begin{Lemma}\label{Lemma1}
For any $\al>-\si$, there is a uniform bound
\[\frac{\GG(z+\si)}{\GG(z)}\ll_{\al,\si}|z+\si|^\si,\quad\Re z\geq\al.\]
\end{Lemma}
\noindent \emph{Proof of Lemma~\ref{Lemma1}.} The function
$\GG(z+\si)/\GG(z)$ is holomorphic in a neighborhood of $\Re
z\geq\al$. For $|z|>2|\si|$ we get, using Stirling's formula,
\[\frac{\GG(z+\si)}{\GG(z)}
\ll_\si\left|\frac{(z+\si)^{z+\si-1/2}}{z^{z-1/2}}\right|
\ll_\si|z+\si|^\si.\] The rest of the values of $z$ (i.e. those
with $\Re z\geq\al$ and $|z|\leq 2|\si|$) form a compact set, so
for these we simply have
\[\frac{\GG(z+\si)}{\GG(z)}\ll_{\al,\si}1\ll_{\al,\si}|z+\si|^\si.\qed\]

\noindent \emph{Proof of Lemma~\ref{Lemma2}.} For $\si=0$ the
statement trivially follows from the definitions (\ref{eq7}) and
(\ref{eq11}). So let $s=\si+it$ where $\si\neq 0$. Pick any
$j=1,\dots,md$ and apply Lemma~\ref{Lemma1} with
\[\al=\frac{1}{2(m^2+1)}-\frac{\si}{2},\quad
z=\frac{1}{4}+\frac{\mu_j}{2}-\frac{\si}{2}+\frac{it}{2}\]
to yield
\[\frac{\GG(1/4+\mu_j/2+\si/2+it/2)}{\GG(1/4+\mu_j/2-\si/2+it/2)}
\ll_{\si,m}|1/4+\mu_j/2+\si/2+it/2|^\si.\]
This is the same as
\[\frac{\GG(1/4+\mu_j/2+s/2)}{\GG(1/4+\overline{\mu}_j/2-s/2)}
\ll_{\si,m}|1/4+\mu_j/2+s/2|^\si.\] Observe that by $|s|\geq\si$ and
$|1/4+\mu_j/2|\gg_m 1$ (see (\ref{eq21})) we have, on the right hand side,
\[|1/4+\mu_j/2+s/2|\ll_{\si,m}|1/4+\mu_j/2||s|.\]
Therefore by taking a product over all $j=1,\dots,md$ we get, using
(\ref{eq7}) and (\ref{eq23}),
\[\left|\pi^{mds}\frac{L(1/2+s,\pi_\infty)}{L(1/2-s,\tilde\pi_\infty)}\right|^{1/2}
\ll_{\si,m,d}\left(\frac{\pi^{md}
C}{N}\right)^{\si/2}|s|^{md\si/2},\quad\Re s=\si.\] By
(\ref{eq11}), this is equivalent to (\ref{eq22}), completing the
proof of Lemma~\ref{Lemma2} and Theorem~\ref{Th1}.\qed

\section{Proof of Theorem~\ref{Th2}}

We can assume that $H(s)$ is the Mellin transform of
$h(x)=-xg'(x)$. Indeed, $h:(0,\infty)\to\RR$ is a smooth function
with functional equation $h(1/x)=h(x)$ and derivatives decaying
faster than any negative power of $x$ as $x\to\infty$, therefore
$H(s)$ is entire and satisfies (\ref{eq20}) and (\ref{eq9}). Also,
\[H(0)=-\int_0^\infty g'(x)=g(0+)=1.\]
Equivalently, $H(s)/s$ is the Mellin transform of $g(x)$, because
by partial integration it follows that
\[-\int_0^\infty g'(x)x^sdx=s\int_0^\infty g(x)x^s\frac{dx}{x}.\]
In any case, $g(x)$ can be expressed as an inverse Mellin transform
\[g(x)=\frac{1}{2\pi i}\int_{(\si)}x^{-s}H(s)\frac{ds}{s}.\]

The idea is to compare $g(x)$ with the function $f(x)$ given by
(\ref{eq25}). We have, for any $\si>0$,
\[f(x)-g(x)=\frac{1}{2\pi i}\int_{(\si)}x^{-s}
\bigl\{C^{-s/2}F(s,\pi_\infty)-1\bigr\}H(s)\frac{ds}{s}.\]
In fact, the integrand is holomorphic in
the entire half plane $\Re s>-1/(m^2+1)$,
so the line of integration can be shifted to any
$\si>-1/(m^2+1)$ without changing the value of the integral.
In particular, the choice $\si=0$ is permissible, i.e.,
\begin{equation}\label{eq24}
f(x)-g(x)=\frac{1}{2\pi i}\int_{-\infty}^\infty x^{-it}
\bigl\{C^{-it/2}F(it,\pi_\infty)-1\bigr\}H(it)\frac{dt}{t}.
\end{equation}

Note that $x^{-it}$ and $C^{-it/2}F(it,\pi_\infty)$ are of modulus
1. For any $\ep>0$, the values of $t$ with
$|t|\geq\min(\eta,C^\ep)$ contribute $O_{\ep,g,m,d}(\eta^{-1})$ to
the integral. This follows from (\ref{eq20}) and $\eta\ll
C^{1/md}$. We shall estimate the remaining contribution via

\begin{Lemma}\label{Lemma3} For any $\ep>0$, there is a uniform bound
\[C^{-it/2}F(it,\pi_\infty)-1\ll_\ep|t|\eta^{-1}C^\ep
,\quad |t|<\min(\eta,C^\ep).\] The implied constant depends only
on $\ep$, $m$ and $d$.
\end{Lemma}

\noindent \emph{Proof.}
As $C^{-it/2}F(it,\pi_\infty)$ lies on the unit circle
it suffices to show that
\[\log\bigl\{C^{-it/2}F(it,\pi_\infty)\bigr\}
\ll_{\ep,m,d}|t|\eta^{-1}C^\ep ,\quad |t|<\min(\eta,C^\ep).\] Here
the left hand side is understood as a continuous function defined
via the principal branch of the logarithm near $t=0$. Using
(\ref{eq7}), (\ref{eq23}) and (\ref{eq11}) we can see that the
derivative (with respect to $t$) of the left hand side is given by
\[\frac{i}{2}\Re\sum_{j=1}^{md}\left\{\frac{\GG'}{\GG}
\left(\frac{1}{4}+\frac{\mu_j}{2}+\frac{it}{2}\right)
-\log\left(\frac{1}{4}+\frac{\mu_j}{2}\right)\right\},\] so we can
further reduce the lemma to
\begin{equation}\label{eq18}
\frac{\GG'}{\GG}\left(\frac{1}{4}+\frac{\mu_j}{2}+\frac{it}{2}\right)
-\log\left(\frac{1}{4}+\frac{\mu_j}{2}\right)=O_{\ep,m,d}(\eta^{-1}C^\ep)
,\quad |t|<\min(\eta,C^\ep).
\end{equation}
Here $1/4+\mu_j/2+it/2$ has real part $\gg_m 1$ by (\ref{eq21})
and absolute value at least $\eta/2$ by (\ref{eq19}).
Therefore a standard bound yields
\[\frac{\GG'}{\GG}\left(\frac{1}{4}+\frac{\mu_j}{2}+\frac{it}{2}\right)=
\log\left(\frac{1}{4}+\frac{\mu_j}{2}+\frac{it}{2}\right)+O_m(\eta^{-1}).\]
For $|t|<\min(\eta,C^\ep)$ we can also see that
\[\log\left(\frac{1}{4}+\frac{\mu_j}{2}+\frac{it}{2}\right)=
\log\left(\frac{1}{4}+\frac{\mu_j}{2}\right)+O(\eta^{-1}C^\ep).\]
It follows from (\ref{eq21}) that $C\gg_{m,d}1$, therefore the
last two estimates add up to (\ref{eq18}) as required.\qed

\medskip

Returning to the integral (\ref{eq24}) it follows from
Lemma~\ref{Lemma3} that the values of $t$ with
$|t|<\min(\eta,C^\ep)$ contribute at most
$O_{\ep,g,m,d}(\eta^{-1}C^{2\ep})$. Altogether we have, by
$C\gg_{m,d}1$,
\[f(x)-g(x)=O_{\ep,g,m,d}(\eta^{-1}C^{2\ep}).\]

We conclude Theorem~\ref{Th2} by combining this estimate with
Corollary~\ref{Cor1} and Molteni's bound (\ref{eq17}).

\end{document}